\documentclass[12pt]{amsart}

\newtheorem{theorem}{Theorem}
\newtheorem{lemma}[theorem]{Lemma}
\newtheorem{proposition}[theorem]{Proposition}
\newtheorem{corollary}[theorem]{Corollary}

\theoremstyle{definition}

\newtheorem{remark}[theorem]{Remark}

\usepackage{amscd,amssymb}
\begin{document}
\baselineskip 15 pt

\title[Automorphic equivalence problem for free algebras]
{Automorphic equivalence problem for free associative algebras of rank two}

\author[Vesselin Drensky and Jie-Tai Yu]
{Vesselin Drensky and Jie-Tai Yu}
\address{Institute of Mathematics and Informatics,
Bulgarian Academy of Sciences, 1113 Sofia, Bulgaria}
\email{drensky@math.bas.bg}
\address{Department of Mathematics, The University of Hong Kong,
Hong Kong SAR, China}
\email{yujt@hkucc.hku.hk}

\thanks
{The research of Vesselin Drensky was partially supported by Grant
MI-1503/2005 of the Bulgarian National Science Fund.}
\thanks{The research of Jie-Tai Yu was partially
supported by a Hong Kong RGC-CERG Grant.}


\begin{abstract}
Let $K\langle x,y\rangle$ be the free associative algebra of rank 2
over an algebraically closed constructive field of any
characteristic. We present an algorithm which decides whether or not
two elements in $K\langle x,y\rangle$ are equivalent under an
automorphism of $K\langle x,y\rangle$. A modification of our
algorithm solves the problem whether or not an element
in $K\langle x,y\rangle$ is a semiinvariant of a nontrivial automorphism.
In particular, it determines whether or not the element has a nontrivial
stabilizer in $\text{Aut}K\langle x,y\rangle$.

An algorithm for
equivalence of polynomials under automorphisms of
${\mathbb C}[x,y]$ was presented by Wightwick. Another, much simpler
algorithm for automorphic equivalence of two polynomials in $K[x,y]$ for any
algebraically closed constructive field $K$ was given
by Makar-Limanov, Shpilrain, and Yu. In our approach we combine
an idea of the latter three authors with an idea
from the unpubished thesis of Lane used to describe
automorphisms which stabilize elements of $K\langle x,y\rangle$.
This also allows us to give a simple proof of
the corresponding result for $K[x,y]$ obtained by
Makar-Limanov,  Shpilrain, and Yu.
\end{abstract}
\maketitle

\section{Introduction}

Let $K$ be an arbitrary field of any characteristic and let $K[x,y]$
and $K\langle x,y\rangle$ be, respectively, the polynomial algebra
in two variables and the free unitary associative algebra of rank 2
(or the algebra of polynomials in the noncommuting variables $x$ and
$y$). Two polynomials $u(x,y)$ and $v(x,y)$ from $K[x,y]$ or
$K\langle x,y\rangle$ are automorphically equivalent, if there
exists an automorphism of the corresponding algebra which brings $u$
to $v$. Wightwick \cite{W} has presented an algorithm which decides
whether or not two polynomials in ${\mathbb C}[x,y]$ are
automorphically equivalent. Makar-Limanov, Shpilrain, and Yu
\cite{MSY} have given a much simpler algorithm which works for
$K[x,y]$ for any algebraically closed constructive $K$. Their method
is based on  {\it peak reduction}. See the survey article of
Shpilrain and Yu \cite{SY2} for other applications of the peak
reduction method to problems for commutative algebra. Shpilrain and
Yu \cite{SY1} have settled a special case of the automorphic
equivalence problem for $K\langle x,y\rangle$, namely, the case
where one of the elements is primitive.

It is a classical result of Jung \cite{J} and van der Kulk \cite{K}
that every automorphism of $K[x,y]$ is tame and is a product of
two kind of automorphisms -- affine and triangular.
Even more, $\text{Aut}K[x,y]$ is isomorphic to $A\ast_CB$, the free product
of the subgroup $A$ of affine automorphisms and
the subgroup $B$ of triangular automorphisms amalgamating
their intersection $C$, the subgroup of affine triangular automorphisms.
This implies that every $\varphi\in\text{Aut}K[x,y]$ has a canonical
form $\varphi=\psi_n\cdots\psi_1$, where each
$\psi_i$ is an affine or triangular automorphism, and the length $n$ is invariant
of $\psi$.

Let $\varphi=\psi_1\cdots\psi_n\in\text{Aut}K[x,y]$ bring $u(x,y)$ to $v(x,y)$.
Makar-Limanov, Shpilrain, and Yu \cite{MSY} have studied the behaviour
of the sequence
\[
d_i=\text{max}(\text{deg}_x(\psi_i\cdots\psi_1u),\text{deg}_y(\psi_i\cdots\psi_1u)),
\quad i=0,1,\ldots,n,
\]
where $\text{deg}_x$ and $\text{deg}_y$ denote the degree with respect to $x$ and $y$,
respectively. If, at some step $d_i\leq d_{i+i}>d_{i+2}$ (a peak),
then they replace
$\psi_{i+1}$ with another affine or triangular automorphism $\psi'_{i+1}$ such that
the new maximum $d'_{i+1}$ of the degrees in $x$ and $y$
of $\psi'_{i+1}\psi_i\cdots\psi_1u$
is smaller than $d_{i+1}$. In this way they move the peak to the right.
This procedure gives that $u(x,y)$ and $v(x,y)$ are automorphically equivalent
if and only if there exist two sequences of affine or triangular automorphisms,
$\rho_1,\ldots,\rho_r$ and $\sigma_1,\ldots,\sigma_s$, with the following property.
The sequences of degrees $p_i=\text{max}(\text{deg}_x(\rho_i\cdots\rho_1u),
\text{deg}_y(\rho_i\cdots\rho_1u))$, $i=1,\ldots,r$, and
$q_j=\text{max}(\text{deg}_x(\sigma_j\cdots\sigma_1v),
\text{deg}_y(\sigma_j\cdots\sigma_1v))$, $j=1,\ldots,s$,
strictly decrease, $p_r=q_s$, and there is an affine automorphism which
sends $\rho_r\cdots\rho_1u$ to $\sigma_s\cdots\sigma_1v$.
The procedure which decides whether or not such sequences of automorphisms exist
reduces the problem to the decision whether or not a system of algebraic equations in
several variables is consistent. Over an algebraically closed constructive $K$
this problem can be solved using Gr\"obner bases techniques.

The $K\langle x,y\rangle$-analogue of the theorem of Jung-van der Kulk
has been established by Czerniakiewicz \cite{Cz} and Makar-Limanov \cite{M}.
Again, every automorphism is tame and $\text{Aut}K\langle x,y\rangle$ is the
free product with amalgamation of the subgroups of triangular and affine
automorphisms. Clearly, the automorphisms of $K\langle x,y\rangle$
fix, up to a nonzero multiplicative constant, the commutator
$[x,y]=xy-yx$.

A theorem of Lane from his unpublished thesis \cite{L} in 1976 states that
an automorphism $\varphi$ of $K\langle x,y\rangle$ has a nontrivial semiinvariant
(i.e., $\varphi u=\lambda u$ for some
$u(x,y)\in K\langle x,y\rangle\backslash \text{span}([x,y]^k\mid k\geq 0)$
and a nonzero constant $\lambda\in K$)
if and only if $\varphi$ is conjugate in $\text{Aut}K\langle x,y\rangle$
to a linear or triangular automorphism.
See Section 9 of Chapter 6 from the book by Cohn \cite{C}
for the improved exposition of the results of Lane.
The idea of the proof is the following.
Every $\varphi=\psi_n\cdots\psi_1\in\text{Aut}K\langle x,y\rangle$
is written in a canonical form and the considerations are
modulo the subspace spanned
by the powers of the commutator $[x,y]$. The first step
is to show that the consecutive action
of nonaffine triangular automorphisms
$\psi_i$ first strictly decrease the total degree of the element
$u(x,y)$. Then, maybe after one action, when the degree is the same,
it starts to increase strictly.
This allows to bound from above
the length $n$ in the canonical form of the automorphisms with $u(x,y)$ as a semiinvariant.
Then the proof is completed by arguments from the theory of free products of groups
with amalgamation.

Lane \cite{L} (see Exercise 6.9.3, p. 362 of \cite{C})
has proved also that the only automorphisms of ${\mathbb C}[x,y]$
with semiinvariants $u(x,y)\in {\mathbb C}[x,y]\backslash {\mathbb C}$ are conjugate
to linear and triangular automorphisms. Smith \cite{S} has determined the eigenvalues
and the eigenvectors of triangular automorphisms. (Clearly, after a linear transformation
of $x$ and $y$, the linear automorphisms also become triangular.)
Recently, the theorem of Lane has been generalized to any field $K$ by
Makar-Limanov, Shpilrain, and Yu \cite{MSY}, involving algebraic geometry.

In the present paper, by combining the algorithmic approach of
Makar-Limanov, Shpilrain, and Yu \cite{MSY} to the automorphic equivalence
in $K[x,y]$ with the idea of Lane (as stated in \cite{C}) in the description
of automorphisms of $K\langle x,y\rangle$ possessing nontrivial semiinvariants,
we obtain an algorithm deciding whether or not two elements
in $K\langle x,y\rangle$ are equivalent under an automorphism
of $K\langle x,y\rangle$. A modification of our algorithm
solves the problem whether or not an element in $K\langle x,y\rangle$
is a semiinvariant of a nontrivial automorphism. In particular,
it determines whether or not the element has a nontrivial stabilizer
in $\text{Aut}K\langle x,y\rangle$.

Our approach works also in the commutative case.
We slightly improve the automorphic equivalence algorithm
of Makar-Limanov, Shpilrain, and Yu \cite{MSY}, replacing
the study of the behaviour of the degree with respect to $x$ and $y$ with
that of the total degree. We simplify also the proof
(over an arbitrary field $K$) of the result
for the stabilizer of $u(x,y)\in K[x,y]\backslash K$,
avoiding usage of algebraic geometry, and
provide an algorithm for the existence of
a nontrivial stabilizer for a given $u(x,y)$.

\section{Preliminaries}

Since automorphisms of $K[x,y]$ and $K\langle x,y\rangle$ are determined by
the images of $x$ and $y$, we shall denote them as
$\varphi=(a,b)$, where $\varphi x=a(x,y)$, $\varphi y=b(x,y)$.
If $\psi=(c,d)$ is another automorphism, we denote their composition as
\[
\psi\varphi=(c,d)(a,b)=(a(c,d),b(c,d)).
\]
The automorphism $\psi$ is affine, if it is of the form
\[
\psi=(\alpha x+\gamma y+\xi,\beta x+\delta y+\eta),\quad
\alpha,\beta,\gamma,\delta,\xi,\eta\in K.
\]
It is triangular, if
\[
\psi=(\alpha x+p(y),\beta y+\eta),\quad
\alpha,\beta\in K^{\ast}=K\backslash 0,\quad \eta\in K,
\]
and the polynomial $p(y)$ does not depend on $x$.
We denote by $A$ and $B$, respectively, the groups of affine and triangular
automorphisms, and with $C=A\cap B$ their intersection. The results
of Jung \cite{J}, van der Kulk \cite{K}, Czerniakiewicz \cite{Cz}, and
Makar-Limanov \cite{M} give that
\[
\text{Aut}K[x,y]\cong \text{Aut}K\langle x,y\rangle\cong A\ast_CB.
\]
Hence $\varphi\in\text{Aut}K[x,y]$
(and similarly for $\varphi\in\text{Aut}K\langle x,y\rangle$)
has the form
\begin{equation}\label{canonical form of automorphisms}
\varphi=\psi_n\cdots \psi_1,
\end{equation}
where each $\psi_i$ is affine or triangular. If two consequent
$\psi_i,\psi_{i+1}$ belong to the same $A$ or $B$, we can replace them
with their product. We may always assume that
if $n>1$ in (\ref{canonical form of automorphisms}),
then either $\psi_i\in A\backslash B$ and $\psi_{i+1}\in B\backslash A$,
or vise versa. We call this decomposition a canonical form of $\varphi$.
The group theoretic properties of $A\ast_CB$ imply that
if $n>1$, then $\varphi\not=1$.
From now on we fix the automorphism
\begin{equation}\label{definition of tau}
\tau=(y,x).
\end{equation}
Then the form (\ref{canonical form of automorphisms})
of the automorphism $\varphi$ can be replaced by
\begin{equation}\label{simplified canonical form of automorphisms}
\varphi=\rho_n\tau\cdots\tau\rho_1\tau\rho_0,
\end{equation}
where $\rho_0,\rho_1,\ldots,\rho_n\in B$ and only $\rho_0$ and $\rho_n$
are allowed to belong to $A$, see for example p. 350 in \cite{C}.
Using the equalities for compositions of automorphisms
\[
(\alpha x+p(y),\beta y+\gamma)=
(x+\alpha^{-1}(p(x)-p(0)),y)(\alpha x+p(0),\beta y+\gamma),
\quad \gamma\in K,
\]
\[
(\alpha x+\xi,\beta y+\eta)\tau=(\beta y+\eta,\alpha x+\xi),\quad \xi,\eta\in K,
\]
we can do further simplifications in (\ref{simplified canonical form of automorphisms}),
assuming that $\rho_i=(x+p_i(x),y)$ with $p_i(0)=0$ for all $i=1,\ldots,n$.

In the next considerations we work in the free algebra $K\langle x,y\rangle$.
We denote by $\vert u(x,y)\vert$ the homogeneous component of maximum
total degree of the nonzero element $u(x,y)\in K\langle x,y\rangle$.
Following Cohn \cite{C}, p. 357,
we call $u(x,y)$ biased if
$\text{deg}_x\vert u\vert\geq \text{deg}_y\vert u\vert$.

Let $V=\text{span}([x,y]^k\mid k\geq 0)$ be the subspace of $K\langle x,y\rangle$
spanned by all powers of the commutator $[x,y]$.
Since $\text{Aut}K\langle x,y\rangle(V)=V$, the group
$\text{Aut}K\langle x,y\rangle$ acts on
the factor vector space $\overline{K\langle x,y\rangle}=K\langle x,y\rangle/V$.
Since $V$ is also graded, $\overline{K\langle x,y\rangle}$
inherits the grading of $K\langle x,y\rangle$. Hence for the nonzero element
$\overline{u(x,y)}\in \overline{K\langle x,y\rangle}$ we may define
$\text{deg }\overline{u}$, $\text{deg}_x\overline{u}$, $\text{deg}_y\overline{u}$,
and $\vert\overline{u}\vert$. Again, $\overline{u(x,y)}$ is biased if
$\text{deg}_x\vert \overline{u}\vert\geq \text{deg}_y\vert \overline{u}\vert$.

The following result is a corollary of a lemma of Lane.

\begin{proposition}\label{lemma of Lane}
{\rm (Corollary 9.6, pp. 361-362 in \cite{C})}
Let $\overline{0}\not=\overline{u(x,y)}\in \overline{K\langle x,y\rangle}$
and let $\rho=(\alpha x+p(y),\beta y+\gamma)$ be a nonaffine triangular automorphism
of $K\langle x,y\rangle$. Then each of the following statements implies the next:

{\rm (i)} $\overline{u(x,y)}$ is biased;

{\rm (ii)} $\text{\rm deg }\overline{u}<\text{\rm deg}(\overline{\tau\rho u})$;

{\rm (iii)} $\text{\rm deg }\overline{u}\leq\text{\rm deg}(\overline{\tau\rho u})$;

{\rm (iv)} $\overline{\tau\rho u} =\overline{u(\alpha y+p(x),\beta x+\gamma)}$
is biased.
\end{proposition}

The following consequence of the proposition is the main step of the proof of
Theorem 6.9.7, p. 361 in \cite{C}. We include the proof for convenience.

\begin{corollary}\label{corollary of lemma of Lane}
Let $\varphi=\rho_n\tau\cdots\tau\rho_1\tau\rho_0\in
\text{\rm Aut}K\langle x,y\rangle$ be written in the form
(\ref{simplified canonical form of automorphisms}).
Let $u(x,y)\in K\langle x,y\rangle\backslash V$ and let
\[
d_{-1}=\text{\rm deg }\overline{u},\quad
d_n=\text{\rm deg}(\overline{\varphi u}),
\]
\[
d_j=\text{\rm deg}(\overline{\tau\rho_j\tau\cdots\tau\rho_1\tau\rho_0u}),\quad
j=0,1,\ldots,n-1.
\]
If $\rho_i,\ldots,\rho_k$ are all nonaffine automorphisms in the decomposition
(\ref{simplified canonical form of automorphisms}), $i\leq 1$, $k\geq n-1$,
then there exists an integer $m$ between $i$ and $k$ such that
\[
d_{-1}=d_{i-1}>d_i>\cdots>d_m\leq d_{m+1}<\cdots <d_k=d_n.
\]
\end{corollary}

\begin{proof}
Clearly, affine automorphisms preserve the degree in $\overline{K\langle x,y\rangle}$.
If $\rho_0$ is affine, then $i=1$ and
$d_{-1}=\text{deg }\overline{u}
=\text{deg}(\overline{\rho_0u})=\text{deg}(\overline{\tau\rho_0u})$. Similarly we conclude
that $d_k=d_n$. Let $m\geq i$ be the largest integer such that
$d_i>d_{i+1}>\cdots>d_m$. Hence either $m=k$ or $m<k$ and $d_m\leq d_{m+1}$.
Applying consecutively parts $\text{(iii)}\implies\text{(iv)}$
and $\text{(i)}\implies\text{(ii)}$ of Proposition \ref{lemma of Lane}
we obtain that $\overline{\tau\rho_{m+1}\tau\cdots\tau\rho_1\tau\rho_0u}$
is biased and $d_{m+1}<d_{m+2}$. We complete the proof by obvious induction.
\end{proof}

We shall also need the following well known lemma, see for example
Lemma 5.1 in \cite{DY}.

\begin{lemma}\label{basis of the free algebra}
As a vector space $K\langle x,y\rangle$ has a basis consisting of the elements
\begin{equation}\label{normal form of noncommutative polynomials}
u_{ab}
=x^{a_1}y^{b_1}[x,y]x^{a_2}y^{b_2}\cdots x^{a_r}y^{b_r}[x,y]x^{a_{r+1}}y^{b_{r+1}},
\end{equation}
where $a_i,b_i,r\geq 0$.
\end{lemma}

Note, that the coefficients of $u(x,y)\in K\langle x,y\rangle$
with resepct to the basis (\ref{normal form of noncommutative polynomials})
can be found explicitly using the equation
$yx=xy-[x,y]$, see e.g. the proof of the lemma in \cite{DY} for details.

\begin{corollary}\label{bounds for the degree}
Let the element $u(x,y)$ in $K\langle x,y\rangle\backslash V$
be written as a linear combination
\[
u(x,y)=\sum\gamma_{ab}u_{ab},
\quad \gamma_{ab}\in K,
\]
of the basis (\ref{normal form of noncommutative polynomials}) and let
$\rho=(\alpha x+p(y),\beta y+\gamma)$ be a nonaffine triangular automorphism of
$K\langle x,y\rangle$. If $a_1=\cdots=a_{q+1}=0$ for all summands $u_{ab}$
with nonzero coefficients $\gamma_{ab}$, then
$\text{\rm deg }u=\text{\rm deg}(\rho u)$. If
some $a_i$ is not equal to $0$ and $\text{\rm deg }p(y)=k>\text{\rm deg }u$,
then $\text{\rm deg}(\overline{\rho u})\geq k$.
\end{corollary}

\begin{proof}
We use the idea of the proof of Theorem 5.2 in \cite{DY}.
If all $a_i$ are equal to 0, then $x$ participates in $u(x,y)$
in commutators $[x,y]$ only. Since $[\alpha x+p(y),\beta y]=\alpha\beta[x,y]$, we obtain
\[
\rho u=\sum\gamma_{ab}\alpha^r\beta^r
(\beta y+\gamma)^{b_1}[x,y](\beta y+\gamma)^{b_2}[x,y]\cdots
[x,y](\beta y+\gamma)^{b_{r+1}}.
\]
Hence $\text{\rm deg }u=\text{\rm deg}(\rho u)$. Now, let some $a_i$ be not equal to 0.
Let $p(y)=\delta_0y^k+\cdots+\delta_{k-1}y+\delta_k$, $\delta_i\in K$, $\delta_0\not=0$,
and $\text{deg }p=k>\text{deg }u$.
We order the elements $u_{ab}$ from the basis
(\ref{normal form of noncommutative polynomials}) lexicographically
assuming that $y>x$. The leading monomial of
$\rho u_{ab}=u_{ab}(\alpha x+p(y),\beta y+\gamma)$ is
\[
(-1)^r\alpha^r\beta^{B+r}\gamma^{kA}
y^{ka_1+b_1+1}xy^{ka_2+b_2+1}x\cdots y^{a_r+b_r+1}xy^{a_{r+1}+b_{r+1}+1},
\]
$A=\sum a_i$, $B=\sum b_i$. Since $k>\text{deg }u\geq A+B+2r$, we obtain that
the different $\rho u_{ab}$ have linearly independent leading monomials.
If $a_i>0$ for some $i$, then the leading monomial of $\rho u_{ab}$ has different
degrees with respect to $x$ and $y$. Hence, the corresponding
bihomogeneous (homogeneous in $x$ and in $y$) component does not belong to $V$.
Since $\text{deg}(\rho u_{ab})\geq kA+B+2r$ and there exists a nonzero $a_i$, we
conclude that $\text{deg}(\overline{\rho u})\geq k$.
\end{proof}

Finally, we need some facts from the theory of Gr\"obner bases.

\begin{proposition}\label{Groebner bases}
Let $K$ be an algebraically closed constructive field and let
$f_j(t_1,\ldots,t_N)$, $j=0,1,\ldots,M$, be a finite set of polynomials
in $K[t_1,\ldots,t_N]$. There is an algorithm which decides whether or not
the system
\[
f_j(t_1,\ldots,t_N)=0, \quad j=1,\ldots,M,
\]
has a solution $(\xi_1,\ldots,\xi_N)\in K^N$
such that $f_0(\xi_1,\ldots,\xi_N)\not=0$.
\end{proposition}

\begin{proof}
The Hilbert Nullstellensatz gives that the system
\[
f_j(t_1,\ldots,t_N)=0,\quad j=1,\ldots,M,
\]
has a solution if and only if
the ideal $I$ of $K[t_1,\ldots,t_N]$ generated by $f_j(t_1,\ldots,t_N)$, $j=1,\ldots,M$,
does not coincide with the whole $K[t_1,\ldots,t_N]$. We can decide
whether or not $I=K[t_1,\ldots,t_N]$ calculating its Gr\"obner basis.
If for every solution $(\xi_1,\ldots,\xi_N)$ of the system we have
$f_0(\xi_1,\ldots,\xi_N)=0$, then the Hilbert Nullstellensatz again implies that
some power of $f_0(t_1,\ldots,t_N)$ belongs to $I$ and $f_0$ belongs to the radical
$\text{Rad}(I)$ of $I$. There is an algorithm which uses Gr\"obner bases and decides whether
or not $f_0(t_1,\ldots,t_N)\in \text{Rad}(I)$, see for example \cite{AL} or the algorithm
RADICALMEMTEST, p. 268 in \cite{BW}.
\end{proof}

\section{The main results}

The following two theorems are the main results of this paper.

\begin{theorem}\label{theorem for equivalence}
Let $K$ be an algebraically closed constructive field
and let $u(x,y),v(x,y)\in K\langle x,y\rangle$. Then there is an algorithm which
decides whether or not $v=\varphi u$ for some
$\varphi=(f(x,y),g(x,y))\in\text{\rm Aut}K\langle x,y\rangle$.
The elements $f(x,y)$ and $g(x,y)$ which determine $\varphi$
can be expressed in terms of the solutions of systems of algebraic equations.
\end{theorem}

\begin{proof}
We want to find $\varphi=(f,g)\in \text{Aut}K\langle x,y\rangle$ such that
$v=\varphi u$.
We can decide efficiently,
presenting $u$ and $v$ as linear combinations of the basis
(\ref{normal form of noncommutative polynomials}) in
Lemma \ref{basis of the free algebra}, whether or not $u(x,y),v(x,y)\in V$.

{\it Case 1.} If
\[
u(x,y)=\sum_{k=0}^m\lambda_k[x,y]^k\in V,\quad \lambda_k\in K,
\]
then $v=\varphi u$ is impossible if $v\not\in V$. Let
\[
v(x,y)=\sum_{k=0}^m\mu_k[x,y]^k\in V,\quad \mu_k\in K.
\]
Since $\varphi [x,y]=\omega[x,y]$, $\omega\in K^{\ast}$,
the action of $\varphi$ on $u$ is determined by the linear
components $f_1=\xi_1x+\xi_2y$ and $g_1=\eta_1x+\eta_2y$
of $f$ and $g$, respectively. Hence
\[
\varphi u=\sum_{k=0}^m\lambda_k\vartheta^k[x,y]^k,\quad
\vartheta=\xi_1\eta_2-\xi_2\eta_1.
\]
Therefore, we have to decide whether or not the equations
\[
t_k(\omega)=\lambda_k\omega^k-\mu_k=0,\quad k=0,1,\ldots,m,
\]
have a common solution. This can be handled efficiently, determining
with the Euclidean algorithm
the greatest common divisor of the polynomials $t_k(\omega)$. It is easy to see that
automorphisms $\varphi$ which send $u$ to $v$ can be characterized
in their normal form (\ref{simplified canonical form of automorphisms})
as follows.
For any common solution $\omega_0$ of the equations $t_k(\omega)=0$ and any
$n\geq 0$ we define in an arbitrary way
$\rho_i=(x+p_i(y),y)$, $p_i(0)=0$, $i=1,\ldots,n$.
Then we choose $\alpha\in K^{\ast}$, $p(y)\in K[y]$, $\gamma\in K$, and define
$\rho_0=(\alpha x+p(y),\alpha^{-1}\omega_0 y+\gamma)$.

{\it Case 2.} Now we assume that $u,v\not\in V$. We repeat the main idea of the
proof of Makar-Limanov, Shpilrain, and Yu \cite{MSY} of the result in the commutative
case. We search for $\varphi$ in the form
$\varphi=\rho_n\tau\cdots\tau\rho_1\tau\rho_0$.
We can efficiently present $u(x,y)$ and $v(x,y)$ in the form
$u=u'+u_V$, $v=v'+v_V$, where $u_V,v_V\in V$  are the sums of the bihomogeneous
components of $u$ and $v$ which are equal, up to multiplicative
constants, to powers of the commutator $[x,y]$.
In the notation of Corollary \ref{corollary of lemma of Lane},
we define
\[
d_{-1}=\text{\rm deg }\overline{u},\quad
d_n=\text{\rm deg}(\overline{v}),
\]
\[
d_j=\text{\rm deg}(\overline{\tau\rho_i\tau\cdots\tau\rho_1\tau\rho_0 u}),\quad
j=0,1,\ldots,n-1.
\]
We assume that both $\rho_0$ and $\rho_n$ are affine. The other cases are
similar and also have to be considered. Hence $i=1$, $k=n-1$. Since $d_{-1}$ and $d_n$ are equal to
the degrees of $\overline{u}$ and $\overline{v}$, they are fixed.
Hence there is a finite number of choices for the sequence of positive integers
$d_j$, $j=0,1,\ldots,n$, with the property that
\[
d_{-1}=d_0>d_1>\cdots>d_m\leq d_{m+1}<\cdots <d_{n-1}=d_n.
\]
This also bounds $n$ from above by
$n\leq\text{deg }\overline{u}+\text{deg }\overline{v}$.
We have to consider all possible sequences $\{d_j\}$.
We fix one of them.
We consider the first and the last automorphisms
$\rho_0=(\xi_0x+\xi'_0y+\xi''_0,\eta_0y+\eta'_0)$
and $\rho_n=(\xi_nx+\xi'_ny+\xi''_n,\eta_ny+\eta'_n)$
with unknown coefficients $\xi_j,\eta_j$, and all other automorphisms
$\rho_j=(\xi_jx+p_j(y),\eta_jy+\eta_j')$ with unknown $\xi_j,\eta_j$ and
unknown polynomials $p_j$.

{\bf Part 1,} {\it Step 1.}
If, writing $u(x,y)$ as linear
combination $u=\sum\gamma_{ab}u_{ab}$
of the basis (\ref{normal form of noncommutative polynomials}), we have
$a_r=0$ for all $\gamma_{ab}\not=0$, then $\rho_0u$ shares the same property.
Hence $\overline{\tau\rho_0u}$ is biased and by Proposition \ref{lemma of Lane},
$d_0<d_1$. Hence $m=0$ and we go to the next part of the procedure.
We assume that there exists a nonzero $a_i$.
It is easy to see, that the same holds
for some basis element in the expression of $\rho_0u$. Then Corollary
\ref{bounds for the degree} gives that the degree of $p_1(y)$ is bounded
by the degree of $u(x,y)$. Let $p_1(y)=\omega_{d_0}y^{d_0}+\cdots+\omega_1y+\omega_0$,
where $\omega_0,\omega_1,\ldots,\omega_{d_0}$ are unknown coefficients.
This bounds the degree
of $\rho_1\tau\rho_0u$ from above in terms of $d_0$,
e.g. $\text{deg}(\rho_1\tau\rho_0u)\leq d_0^2$.
We write $\rho_1\tau\rho_0u$ in the form
\[
\rho_1\tau\rho_0u=\sum \delta_i z_{i_1}z_{i_2}\cdots z_{i_s},\quad
z_{i_j}=x,y,
\]
where the coefficients $\delta_i=\delta(\xi,\eta,\omega)$
are polynomials in $\xi_j,\eta_j,\omega_j$.
Now we use the equalities
$\text{deg}(\overline{\rho_1\tau\rho_0u})=d_1$ and
$\text{deg}(\rho_1\tau\rho_0u)\leq d_0^2$.
The monomials $z_{i_1}z_{i_2}\cdots z_{i_s}\in K\langle x,y\rangle$ of
degree $s>d_0^2$ do not participate in $\rho_1\tau\rho_0u$.
If $d_1<s\leq d_0^2$ and $\text{deg}_x(z_{i_1}z_{i_2}\cdots z_{i_s})$ is different from
$\text{deg}_y(z_{i_1}z_{i_2}\cdots z_{i_s})$, then $\delta_i=0$.
If $\text{deg}_x(z_{i_1}z_{i_2}\cdots z_{i_s})
=\text{deg}_x(z_{i_1}z_{i_2}\cdots z_{i_s})$ for $d_1<s\leq d_0^2$, then we write
the corresponding bihomogeneous component in the form
$\sum\delta_{i'}z_{i_1}z_{i_2}\cdots z_{i_s}=\vartheta [x,y]^{s/2}$ with
unknown coefficient $\vartheta$ and again obtain equations of the form
$\delta_i=0$ or $\delta_i-\pm\vartheta$. In this way we obtain a finite system
of algebraic equations
\begin{equation}\label{system to solve}
\Delta_q(\xi,\eta,\omega,\vartheta)=0,\quad q=1,\ldots,Q.
\end{equation}
We want to decide whether or not the system has a solution with the property
that $\text{deg}(\overline{\rho_1\tau\rho_0u})=d_1$, the coefficients
$\xi_0,\eta_0,\xi_1,\eta_1$ are nonzero, and the polynomial $p_1(y)$ is of degree
$\geq 2$. This can be done effectively using Proposition \ref{Groebner bases}.

{\it Step 2.} We repeat Step 1 with the element $\tau\rho_1\tau\rho_0u$ instead of
with $u$. If the element $\overline{\tau\rho_1\tau\rho_0u}$ is biased, we have $m=1$ and go
to the next part of the procedure. If $\overline{\tau\rho_1\tau\rho_0u}$ is not biased,
then we bound from above the degree of the polynomial $p_2(y)$ in the definition
of $\rho_2$ and, continuing as above, add new equations to the system of algebraic
equations (\ref{system to solve}).

We continue till the $m-1$'st step,
and obtain the polynomial
\[
w_1=\rho_{m-1}\tau\cdots\rho_1\tau\rho_0u.
\]
when $d_m\leq d_{m+1}$. We finish this part of the procedure.

{\bf Part 2.} We start
a similar procedure with $v(x,y)$, applying to it
$\rho_j^{-1}\tau\cdots\tau\rho_{n-1}^{-1}\tau\rho_n^{-1}$ for $j=n,n-1,\ldots,m$
if $d_m<d_{m+1}$ and for $j=n,n-1,\ldots,m+1$ if $d_m=d_{m+1}$.

{\bf Part 3.} If $d_m<d_{m+1}$,
we obtain the element
\[
w_2=\rho_m^{-1}\tau\cdots\tau\rho_{n-1}^{-1}\tau\rho_n^{-1}v,
\]
and a system of algebraic equations depending on the unknown coefficients of
$\rho_i$.
Since $w_1=\tau w_2$, we obtain one more relation between the coefficients
of $\rho_0,\rho_1,\ldots,\rho_n$. If $d_m=d_{m+1}$, then we
consider
\[
w_2=\rho_{m+1}^{-1}\tau\cdots\tau\rho_{n-1}^{-1}\tau\rho_n^{-1}v.
\]
Then $w_2=\tau\rho_m\tau w_1$ and we have two possibilities.
If, writing $\tau w_1$ as a linear
combination $\tau w_1=\sum\gamma'_{ab}u_{ab}$
of the basis (\ref{normal form of noncommutative polynomials}), we have
$a_j=0$ for all $\gamma'_{ab}\not=0$, then
\[
\tau w_1=\sum\gamma'_{ab}
y^{b_1}[x,y]\cdots [x,y]y^{b_{r+1}},
\]
\[
\rho_m\tau w_1=\sum\gamma'_{ab}\xi_m^r\eta_m^r
(\eta_my+\eta_m')^{b_1}[x,y]\cdots [x,y](\eta_my+\eta_m')^{b_{r+1}}.
\]
Hence the result does not depend on the polynomial $p_m(y)$,
we can choose it to be arbitrary.
The corresponding algebraic system does not depend on its coefficients.
If some $a_i$ is positive, then we bound the degree of $p_m(y)$ and determine
whether or not the obtained system has a solution with nonzero $\xi_j,\eta_j$
and nonlinear $p_j(y)$.
\end{proof}

\begin{theorem}\label{theorem for conjugation}
Let $K$ be an algebraically closed constructive field
and let $u(x,y)\in K\langle x,y\rangle$. Then there is an algorithm which
decides whether or not $u$ is a semiinvariant of some
$\varphi\in\text{\rm Aut}K\langle x,y\rangle$,
and $\varphi$
can be expressed in terms of the solutions of algebraic systems.
\end{theorem}

\begin{proof} Let $0\not=u(x,y)\in K\langle x,y\rangle$.
We want to find $\varphi\in \text{Aut}K\langle x,y\rangle$ and a constant $\lambda$
such that $\varphi u=\lambda u$.
If $u(x,y)\in V$, then the action of $\varphi=(f,g)$ is determined
by $\vartheta=\xi_1\eta_2-\xi_2\eta_1$, where
$f_1=\xi_1x+\xi_2y$ and $g_1=\eta_1x+\eta_2y$ are the linear components
of $f$ and $g$, respectively.
In particular, $u$ is stabilized by any $\varphi$ with $\vartheta=1$.
We can find all possible values of $\vartheta$ as in the first part of the proof of
Theorem \ref{theorem for equivalence}. If
$\varphi u=\sum_{k=0}^m\lambda_k\vartheta^k[x,y]^k$,
then $\varphi u=\lambda u$ if and only if
$\lambda=\vartheta^k$ for all $k$ such that $\lambda_k\not=0$.
If $u(x,y)\not\in V$, then the theorem of Lane
implies that $\varphi=\psi^{-1}\rho\psi$ for some triangular or
affine $\rho$ and some $\psi\in \text{Aut}K\langle x,y\rangle$. We
write $\psi$ in the form (\ref{simplified canonical form of
automorphisms}), $\psi=\rho_n\tau\cdots\tau\rho_1\tau\rho_0$. If
$\rho$ is triangular, we have
\[
\varphi=(\rho_n\tau\cdots\tau\rho_1\tau\rho_0)^{-1}
\rho(\rho_n\tau\cdots\tau\rho_1\tau\rho_0)
\]
\[
=\rho_0^{-1}\tau\rho_1^{-1}\tau\cdots\tau(\rho_n^{-1}\rho\rho_n)\tau\cdots\rho_1\tau\rho_0.
\]
If $\rho$ is nontriangular affine, then it has the form $\rho=\rho''\tau\rho'$
some some affine triangular $\rho',\rho''$ and we proceed in a similar way.
Then we complete the proof as in the second case of the proof of
Theorem \ref{theorem for equivalence}.
\end{proof}

\begin{remark} The unitarity of $K\langle x,y\rangle$ is not essential.
The same proofs work in the free nonunitary associative algebra in two variables.
\end{remark}

\section{Applications to commutative case}

From now on we work in the polynomial algebra $K[x,y]$
over an arbitrary field $K$,
keeping some notation from the case $K\langle x,y\rangle$.
If $0\not=u(x,y)\in K[x,y]$, we denote by $\vert u\vert$ the
homogeneous component of maximum total degree and say that $u$ is biased
if $\text{deg}_x\vert u\vert\geq\text{deg}_y\vert u\vert$.

One of the main steps in the approach of Makar-Limanov, Shpilrain,
and Yu \cite{MSY} is Lemma 2 in \cite{MSY}. A minor modification in
its proof allows us to simplify the proof of the theorem for the
existence of a nontrivial stabilizer of $u\in K[x,y]$.

\begin{proposition}\label{commutative lemma of Lane}
Let $u(x,y)\in K[x,y]\backslash K$ be biased and let
$\rho=(\alpha x+p(y),\beta y+\gamma)$ be a nonaffine triangular automorphism.
Then $\text{\rm deg}(\rho u)>\text{\rm deg }u$.
\end{proposition}

\begin{proof}
For simplicity of the exposition we assume that $\alpha=\beta=1$ and
$p(x)=x^k+\pi_{k-1}x^{k-1}+\cdots+\pi_1x+\pi_0$ is monic, with $k\geq 2$.
Let the homogeneous component of maximum degree of $u(x,y)$ be
\[
\vert u\vert=\gamma_ax^ay^b+\gamma_{a-1}x^{a-1}y^{b+1}
+\cdots+\gamma_1xy^{a+b-1}+\gamma_{a+b}y^{a+b},\quad \gamma_a\not=0.
\]
Since $u(x,y)$ is biased, we have $a\geq b$. Define a $(k,1)$-grading on
$K[x,y]$ assuming that $\text{deg}_{(k,1)}x=k$, $\text{deg}_{(k,1)}y=1$. Let
$u_m(x,y)$ be the homogeneous component of $u$ of $(k,1)$-degree $m$, and let
\[
u_{ka+b}(x,y)=\beta_dx^dy^j+\beta_{d-1}x^{d-1}y^{j+k}+\cdots+\beta_ex^ey^{j+k(d-e)},
\quad \beta_d,\beta_e\not=0.
\]
Then, over the algebraic closure $\overline{K}$ of $K$,
the $(k,1)$-homogenity of $u_{ka+b}$ implies the decompostion
\[
u_{ka+b}(x,y)=\xi x^qy^r(x-y^k)^s\prod_{i=1}^t(x-\lambda_iy^k),
\quad 1\not=\lambda_i\in\overline{K}.
\]
Clearly,
\[
\text{deg }u_{ka+b}=q+r+k(s+t)=e+(j+k(d-e))
\]
\[
>e+1+(j+k(d-e-1))>\cdots>d+j.
\]
Since $k\geq 2$, the only summand of maximum total degree contained in $u_{ka+b}$
is $\gamma_ax^ay^b$. We conclude that $\gamma_ax^ay^b=\beta_ex^ey^{j+k(d-e)}$
and this implies that
\[
(a,b)=(e,j+k(d-e))=(q,r+k(s+t)),
\]
\[
a+b=\text{deg }u=\text{deg }u_{ka+b}
=q+r+k(s+t).
\]
As in the proof of Lemma 2 \cite{MSY}, the first step is to show that
$\text{deg }u_{ka+b}(x+y^k,y)>\text{deg }u_{ka+b}(x,y)$. Let us assume that the oposite
inequality $\text{deg }u_{ka+b}(x+y^k,y)\leq\text{deg }u_{ka+b}(x,y)$ holds.
Since
\[
u_{ka+b}(x+y^k,y)=\xi (x+y^k)^qy^rx^s\prod_{i=1}^t(x-(\lambda_i-1)y^k),
\]
and $k\geq 2$, we derive that
\[
\text{deg }u_{ka+b}(x+y^k,y)=kq+r+s+kt\leq q+r+k(s+t)=\text{deg }u_{ka+b}(x,y),
\]
\[
(k-1)q\leq (k-1)s,\quad q\leq s.
\]
Since $u(x,y)$ is biased and
\[
q=a\geq b=r+k(s+t),\quad q>0,r,s,t\geq 0,
\]
we obtain that $q\geq 2s$. This is a contradiction because
we already have $s\geq q>0$.
In this way $\text{deg }u_{ka+b}(x+y^k,y)>\text{deg }u(x,y)$.
Since the leading $(k,1)$-components of
$x+p(y)$ and $y+\gamma$ are $x+y^k$ and $y$, respectively, we derive that
\[
u_m(x+p(y),y+\gamma)=u_m(x+y^k,y)
\]
\[
+\text{$(k,1)$-homogeneous components of lower
$(k,1)$-degree}.
\]
Hence the monomials of $u_{ka+b}(x+y^k,y)$ can vanish in
$u(x+p(y),y+\gamma)$ only if they cancel with some monomials
from $u_m(x+p(y),y+\gamma)$ for $m>ka+b$. Let $m_0$ be the $(k,1)$-degree of
$u(x,y)$. If $m_0=ka+b$, then the monomials of $u_{ka+b}(x+y^k,y)$
do not cancel with anything.
Hence $\text{deg}(\rho u)\geq\text{deg }u_{ka+b}(x+y^k,y)>\text{deg }u(x,y)$.
So, we may assume that $m_0>ka+b$.
Again, the leading $(k,1)$-component of $u_{m_0}(x+p(y),y+\gamma)$ is
$u_{m_0}(x+y^k,y)$. If
\[
u_{m_0}(x,y)=\xi_0 x^{q_0}y^{r_0}(x-y^k)^{s_0}\prod_{i=1}^{t_0}(x-\lambda'_iy^k),
\quad \lambda'_i\not=1,
\]
then the leading $(k,1)$-component of $u_{m_0}(x+p(y),y+\gamma)$ is
\[
u_{m_0}(x+y^k,y)=\xi_0 (x+y^k)^{q_0}y^{r_0}x^{s_0}\prod_{i=1}^{t_0}(x-(\lambda'_i-1)y^k)
\]
and does not cancel with other elements of $u(x+p(y),y+\gamma)$.
In particular, $\text{deg }u(x+p(y),y+\gamma)\geq \text{deg }u_{m_0}(x+y^k,y)$.
We have the inequalities
\begin{equation}\label{first inequality}
a+b=\text{deg }u(x,y)\geq \text{deg }u_{m_0}(x,y)=q_0+r_0+k(s_0+t_0),
\end{equation}
\begin{equation}\label{second inequality}
kq_0+r_0+k(s_0+t_0)=m_0>ka+b.
\end{equation}
The sum of (\ref{first inequality}) and (\ref{second inequality}) gives
\[
a+b+kq_0+r_0+k(s_0+t_0)>ka+b+q_0+r_0+k(s_0+t_0),
\]
\[
(k-1)q_0>(k-1)a,\quad q_0>a,
\]
If we assume that
\begin{equation}\label{third inequality}
a+b=\text{deg }u(x,y)\geq\text{deg }u_{m_0}(x+y^k,y)=kq_0+r_0+s_0+kt_0,
\end{equation}
then the sum of (\ref{second inequality}) and (\ref{third inequality}) implies
\[
a+b+kq_0+r_0+k(s_0+t_0)>ka+b+kq_0+r_0+s_0+kt_0,
\]
\[
(k-1)s_0>(k-1)a,\quad s_0>a,
\]
and (\ref{first inequality}) gives
\[
2a\geq a+b\geq q_0+r_0+k(s_0+t_0)\geq q_0+s_0>2a,
\]
which is impossible. Hence
\[
\text{deg }u(x+p(y),y+\gamma) \geq\text{deg }u_{m_0}(x+y^k,y)>\text{deg }u(x,y).
\]
\end{proof}

Proposition \ref{commutative lemma of Lane} implies immediately
commutative analogues of Proposition \ref{lemma of Lane}
and Corollary \ref{corollary of lemma of Lane}. We shall state the first of them.

\begin{corollary}\label{commutative corollary}
Let $u(x,y)\in K[x,y]\backslash K$
and let $\rho=(\alpha x+p(y),\beta y+\gamma)$ be a nonaffine triangular automorphism
of $K[x,y]$. Then each of the following statements implies the next:

{\rm (i)} $u(x,y)$ is biased;

{\rm (ii)} $\text{\rm deg }u<\text{\rm deg}(\tau\rho u)$;

{\rm (iii)} $\text{\rm deg }u\leq\text{\rm deg}(\tau\rho u)$;

{\rm (iv)} $\tau\rho u=u(\alpha y+p(x),\beta x+\gamma)$
is biased.
\end{corollary}

\begin{proof}
The only part of the proof left is the implication
$\text{(iii)}\implies \text{(iv)}$. If some $v(x,y)\in K[x,y]$ is not biased, then
$v(y,x)=\tau v$ is. Hence, if
$\tau\rho u$ is not biased, then
$\rho u$ is and Proposition \ref{commutative lemma of Lane} gives that
$\text{deg }u=\text{deg}(\rho^{-1}(\rho u))>\text{deg}(\rho u)=\text{deg}(\tau\rho u)$
which is a contradiction.
\end{proof}

Now we can prove easily the theorem of Lane \cite{L} and Makar-Limanov, Shpirain, and Yu
\cite{MSY}.

\begin{theorem}
If $K$ is any field and the automorphism $\varphi$ of $K[x,y]$ is not conjugate to a
linear or triangular automorphism, then any semiinvariant
$u(x,y)\in K[x,y]$ of $\varphi$ is a constant.
\end{theorem}

\begin{proof}
Let $u(x,y)\in K[x,y]\backslash K$ and let $G$ be the subgroup of $\text{Aut}K[x,y]$
which stabilizes the vector space spanned by $u(x,y)$. Writing $\varphi\in G$
in the form (\ref{simplified canonical form of automorphisms}),
$\varphi=\rho_n\tau\cdots\tau\rho_1\tau\rho_0$ and applying the commutative analogue of
Corollary \ref{corollary of lemma of Lane}, we obtain that the length $n$ in the
expression of $\varphi$ is bounded by $2\cdot\text{deg }u$. Now the proof is completed
by the well known theorem in group theory (see e.g. Theorem 6.8.7, p. 351 \cite{C}),
which states that if $G$ is a subgroup of $A\ast_CB$ and its elements are of the
form
\[
g=a_mb_m\cdots a_1b_1,\quad a_i\in A,\quad b_i\in B,
\]
where the integers $m=m(g)$ are bounded by the same $n$ for all $g\in G$,
then $G$ is conjugate to a subgroup
of $A$ or $B$. In order to replace the affine automorphisms with linear ones,
we need to use the fact that $\text{Aut}K[x,y]$ is also a free product
of the linear group $GL_2(K)$ and the triangular group $B$
with amalgamation over their intersection.
\end{proof}

Clearly, we have analogues for $K[x,y]$ of the algorithms described
in Theorems \ref{theorem for equivalence} and \ref{theorem for conjugation}
(compare the first algorithm with this of Makar-Limanov, Shpilrain, and Yu \cite{MSY}).
In particular, when $K$ is an algebraically closed constructive field,
we can decide whether or not $u\in K[x,y]$ is a semiinvariant of some
$\varphi\in\text{\rm Aut}K[x,y]$ and to express
$\varphi$ in terms of solutions of algebraic systems.

\begin{remark}
Clearly, over an algebraically closed field $K$
any linear automorphism can be triangularized.
Smith \cite{S} has determined the eigenvalues and the eigenvectors of
any triangular automorphism $\rho$ of $K[x,y]$ when $\text{char}K=0$.
Up to conjugation, the possibilities are:

(i) $\rho=(\alpha x,\beta y)$, $u(x,y)$ is a linear combination of monomials
$x^ny^m$ with the same value of $\alpha^n\beta^m$;

(ii) $\rho=(\alpha x,\beta y+\gamma)$, $\gamma\not=0$, $u(x,y)$ does not depend on
$y$ and is a linear combination of powers $x^n$ with the same value of $\alpha^n$;

(iii) $\rho=(\alpha x+p(y),\beta y)$, $p(y)\not=0$, $u(x,y)$ does not depend
on $x$ and is a linear combination of powers $y^m$ with the same value of $\beta^m$.

If $u=w(f)$ for some coordinate $f(x,y)$ and some polynomial $w(z)$, we have
\[
u_x=\frac{\partial u}{\partial x}=w'(f)f_x,\quad
u_y=\frac{\partial u}{\partial y}=w'(f)f_y,
\]
and the ideal of $K[x,y]$ generated by $f_x$ and $f_y$
conicides with the whole $K[x,y]$. Hence the greatest common divisor
of $u_x$ and $u_y$ is $w'(f)$
and this can be used to determine whether or not $u$ is a semiinvariant
of a nontrivial automorphism in the cases (ii) and (iii). We cannot see how to handle
directly the case (i), i.e., to determine whether or not $u=w(f,g)$ with some specific
properties of the polynomial $w(z,t)$.
\end{remark}

\section*{Acknowledgements}

The  authors would like to thank L. Makar-Limanov and V. Shpilrain
for helpful comments and an anonymous referee for the numerous
suggestions for improving the exposition.

\end{document}